\documentclass[11pt]{amsart}
\usepackage{amsmath,amssymb,amsthm,mathtools}
\usepackage[T1]{fontenc}
\usepackage{lmodern}
\usepackage{microtype}
\usepackage[hidelinks]{hyperref}
\numberwithin{equation}{section}
\newtheorem{theorem}{Theorem}[section]
\newtheorem{proposition}[theorem]{Proposition}
\newtheorem{lemma}[theorem]{Lemma}

\theoremstyle{definition}
\newtheorem{definition}[theorem]{Definition}
\newtheorem{example}[theorem]{Example}
\newtheorem{remark}[theorem]{Remark}
\newcommand{\J}{\mathcal J}
\newcommand{\OO}{\mathcal O}
\newcommand{\aideal}{\mathfrak a}
\newcommand{\supp}{\mathrm{supp}}
\newcommand{\Q}{\mathbb Q}
\newcommand{\Z}{\mathbb Z}
\newcommand{\R}{\mathbb R}
\newcommand{\F}{\mathbb F}
\newcommand{\C}{\mathbb C}
\newcommand{\xx}{\boldsymbol{x}}
\newcommand{\one}{\boldsymbol{1}}

\title[Uniform comparison for affinely independent support]{Uniform comparison of multiplier and test ideals for polynomials with affinely independent support}
\author{Takafumi Shibuta}
\address{Department of Information Science,
Faculty of Science and Engineering,
Kyushu Sangyo University}
\date{}
\thanks{This work has been supported by JSPS Grant-in-Aid for Scientific Research (C) (No. 22K03334).}
\email{tshibuta@ip.kyusan-u.ac.jp}
\subjclass[2020]{Primary 13A35; Secondary 14F18, 14M25}
\keywords{test ideal, multiplier ideal, reduction modulo primes}

\begin{document}
\begin{abstract}
Let $f\in\Q[x_1,\ldots,x_n]$ have affinely independent exponent vectors.
We prove that its test ideals agree with the reductions of its multiplier
ideals at all nonnegative real exponents for every sufficiently large prime in one
arithmetic progression. The progression depends only on the support of $f$.
\end{abstract}
\maketitle

\section{Introduction}

Multiplier ideals in characteristic zero and test ideals in positive
characteristic measure the singularities of ideals. For a nonzero ideal
$I\subseteq\Q[\xx]$, the comparison results of \cite{HY} and \cite[Section~4]{MS} show that, after choosing
a common model, $\tau(I_p^t)\subseteq\J(I^t)_p$ for all sufficiently large
primes $p$ and all real $t\ge0$, where the subscript $p$ denotes reduction
modulo $p$. For each fixed $t$, equality holds for all sufficiently large
$p$, with a bound that may depend on $t$.
The comparison conjecture of Musta\c{t}\u{a} and Srinivas
\cite[Conjecture~1.2]{MS} states that there are infinitely many
primes $p$ for which equality holds simultaneously for all real $t\ge0$. We prove this conjecture
for principal ideals generated by polynomials with affinely independent support.

Write
\[
 f=\sum_{j=1}^s\gamma_j\xx^{u_j}\in\Q[\xx]=\Q[x_1,\dots,x_n],
\]
where $\gamma_j\in\Q^*$ and the $u_j\in\Z_{\ge 0}^n$ are distinct vectors. 
We call $\supp f=\{u_1,\dots ,u_s\}$ the support of $f$. 
Let
$\aideal_f=\langle\xx^{u_1},\ldots,\xx^{u_s}\rangle$ be the term ideal of $f$.
The $u_j$ are affinely independent if and only if $(u_1,1),\dots,(u_s,1)$ are linearly independent over $\Q$. 

\begin{theorem}\label{thm:main}
Suppose that the support of $f$ is affinely independent.
There is a positive integer $N$, depending only on the support, such that
for every sufficiently large prime number $p$ satisfying $p\equiv1\pmod N$, 
\[
 \tau(f_p^{t})=\J(f^{t})_p
   =\J(\aideal_f^t)_p\qquad(0\le t<1).
\]
By Skoda's theorem, the comparison holds for every real exponent $t\ge1$ as well.
Thus the comparison conjecture holds for $\langle f\rangle$.
\end{theorem}

For $0\le t<1$, the general-member formula and a torus action give
$\J(f^t)=\J(\aideal_f^t)$. It therefore suffices to prove
$\J(\aideal_f^t)_p\subseteq\tau(f_p^t)$ for all $0\le t<1$ simultaneously
at each sufficiently large prime in the stated progression.

Musta\c{t}\u{a}, Takagi and Watanabe \cite[Example~4.2]{MTW} compare the
$F$-pure and log canonical thresholds under the same affine-independence
hypothesis. Their proof uses rational points of a polytope to choose the multiplicities
in a multinomial expansion. We use the same idea here.
Hern\'andez \cite{Her} develops splitting polytopes and base-$p$ methods
for threshold comparisons under related hypotheses.
Hern\'andez, Teixeira and Witt \cite[Theorem~6.4]{HTW} prove simultaneous
comparison at all exponents for diagonal polynomials. Their
Proposition~6.2 separates terms in a Frobenius basis; their Remark~6.3
also identifies monomial test ideals of a polynomial with generalized
Frobenius powers of its term ideal.

Hern\'andez, Teixeira and Witt developed the theory of generalized
Frobenius powers. For each fixed $t\in(0,1)$, they prove in
\cite[Theorem~6.3]{HTWF} that these powers agree with the reductions
of multiplier ideals for all sufficiently large primes. In
\cite[Question~6.6]{HTWF}, they ask whether this agreement holds for
all $t\in(0,1)$ simultaneously on an infinite set of primes.

Theorem~\ref{thm:main} extends the all-exponent comparison for diagonal
polynomials to arbitrary affinely independent support. The arithmetic
construction follows \cite[Example~4.2]{MTW}; the additional step is to
choose rational points for finitely many monomial generators at the right
endpoints of the constancy intervals. One congruence condition then works
for every interval.
%%%%%%%%%%%%%%%%%%%%%%%%%%%%%%%%%%%%%%
%%%%%%%%%%%%%%%%%%%%%%%%%%%%%%%%%%%%%%
\section{Preliminaries}
For $a\in\Z_{\ge 0}^n$, write $\xx^a=x_1^{a_1}\cdots x_n^{a_n}$, and put
$\one=(1,\ldots,1)$. For real vectors, floors are taken coordinatewise,
and $a<b$ (respectively $a\le b$) means $a_i<b_i$ (respectively
$a_i\le b_i$) for every $i$.

Let $X=\operatorname{Spec}\Q[\xx]$ and let $I$ be a nonzero ideal of $\Q[\xx]$.
Choose a log resolution $\pi:Y\to X$ of $I$. Thus $Y$ is smooth,
$I\OO_Y=\OO_Y(-F)$, and the support of $F$ together with the exceptional
divisor has simple normal crossings. For real $t\ge0$, the multiplier
ideal is
\[
 \J(I^t)=\pi_*\OO_Y(K_{Y/X}-\lfloor tF\rfloor),
\]
where $K_{Y/X}$ is the relative canonical divisor. We identify this
ideal sheaf with its ideal of global sections. It is independent of the
resolution \cite[Section~3.1]{MS}.

For an ideal $I\subseteq \F_p[\xx]$, let
$I^{[q]}$, $q=p^e$, be generated by the $q$-th powers of its elements, and let
$I^{[1/q]}$ be the smallest ideal $L$ with $I\subseteq L^{[q]}$.

Since $\F_p[\xx]$ is regular and $F$-finite, we use the
Frobenius-root definition of Blickle, Musta\c{t}\u{a} and Smith
\cite[Definition~2.9]{BMS} for generalized test ideals. 
In this setting it agrees with the definition of Hara and Yoshida \cite{HY}, as shown in
\cite[Proposition~2.22]{BMS}.

\begin{definition}[Test ideal]
Let $I\subseteq \F_p[\xx]$ be a nonzero ideal and let $t\ge0$ be real.
The test ideal of $I$ with exponent $t$ is
\begin{equation}\label{eq:test-union}
 \tau(I^t)=\bigcup_{e\ge1}
          \bigl(I^{\lceil tp^e\rceil}\bigr)^{[1/p^e]}.
\end{equation}
\end{definition}
The ideals $\bigl(I^{\lceil tp^e\rceil}\bigr)^{[1/p^e]}$, $e\ge 1$, in \eqref{eq:test-union} form an increasing sequence and hence
stabilize. In particular, each term is contained in $\tau(I^t)$.
Frobenius roots and test ideals preserve inclusions of ideals.

For $\eta_a\in\F_p$ and $h_a\in\F_p[\xx]$, the Frobenius basis formula
\cite[Proposition~2.5]{BMS} gives
\[
 H=\sum_{a\in\{0,\ldots,q-1\}^n}\eta_a h_a^q\xx^a
 \quad\Longrightarrow\quad
 \langle H\rangle^{[1/q]}=\langle h_a\mid\eta_a\ne0\rangle.
\]

\begin{lemma}\label{lem:froot}
Fix $q=p^e$ and write $H=\sum_{j=1}^{\ell}\eta_j\xx^{v_j}$,
where $\eta_j\in\F_p$ and $v_j\in\Z_{\ge 0}^n$.
If $v_i\not\equiv v_j\pmod q$ whenever $i\ne j$, then
\[
 \langle H\rangle^{[1/q]}
 =\langle \xx^{\lfloor v_j/q\rfloor}\mid \eta_j\ne0\rangle.
\]
\end{lemma}

\begin{proof}
Write $v_j=q\lfloor v_j/q\rfloor+a_j$ with $a_j\in\{0,\ldots,q-1\}^n$.
The $a_j$ are distinct, so the preceding formula applies with
$h_{a_j}=\xx^{\lfloor v_j/q\rfloor}$.
\end{proof}

Let
\begin{equation}\label{eq:f}
f=\sum_{j=1}^s\gamma_j\xx^{u_j}\in\Q[\xx].
\end{equation}
Choose $A=\Z[1/D]$ so that every coefficient $\gamma_j$ is a unit in $A$.
For $p\nmid D$, reduction of coefficients gives $f_p\in\F_p[\xx]$.
For $m\in\Z_{\ge 0}$ and $0\le t<1$, the principal case of Skoda's theorem gives
\begin{equation}\label{eq:skoda}
 \J(f^{m+t})=f^m\J(f^t),\qquad
 \tau(f_p^{m+t})=f_p^m\tau(f_p^t).
\end{equation}
The first identity follows by extracting the integral principal divisor (see \cite[Section~11.1.A]{Laz}).
For the second, use $\lceil q(m+t)\rceil=mq+\lceil qt\rceil$ and
$\langle f_p^{mq}H\rangle^{[1/q]}=f_p^m\langle H\rangle^{[1/q]}$ in
\eqref{eq:test-union}; see also \cite[Proposition~2.25]{BMS}.

For an ideal $I$ with model $I_A\subseteq A[\xx]$, write $I_p$ for
its image in $\F_p[\xx]$. We choose the models of the multiplier ideals
using a single log resolution as in \cite[Section~3.2]{MS}.
After enlarging $D$, this gives compatible reductions for all exponents:
there are only finitely many multiplier ideals for $0\le t<1$, and
\eqref{eq:skoda} determines the rest.

\begin{theorem}[Hara--Yoshida comparison, principal case]
For every sufficiently large prime $p$,
$\tau(f_p^t)\subseteq\J(f^t)_p$ for every real $t\ge0$.
For each fixed real $t\ge0$, equality holds for every sufficiently
large prime $p$, with a bound that may depend on $t$.
\end{theorem}
For the comparison theorem see \cite{HY}. The uniform inclusion and the
fixed-exponent statement in this form are given in \cite[Section~4]{MS}.

%\begin{remark}\label{rem:large}
%Suppose the support of $f$ is affinely independent, and let $\Delta>0$ be the
%greatest common divisor of the maximal minors of the matrix with columns
%$(u_j,1)$. If $p\nmid\Delta$, some maximal minor is a unit modulo every
%$q=p^e$. The matrix therefore gives an injection
%$(\Z/q\Z)^s\to(\Z/q\Z)^{n+1}$ for every $e\ge1$.
%Thus $p>\max\{D,\Delta\}$ suffices for reduction with the same support
%and for injectivity. Increasing this bound to include the bound in the
%comparison theorem ensures the inclusion for every $t$ as well.
%\end{remark}

We use Howald's formula for monomial multiplier ideals \cite{How01}
and the corresponding test ideal formula of Hara--Yoshida
\cite[Theorem~6.10(3)]{HY}. The rational-exponent formulas extend to
real exponents by right continuity (for test ideals, see
\cite[Proposition~2.14]{BMS}).

\begin{proposition}\label{prop:monomial}
For $\aideal=\langle\xx^{u_1},\ldots,\xx^{u_s}\rangle$ and $t\ge0$,
\begin{equation}\label{eq:monomial}
 \J(\aideal^t)=
 \left\langle\xx^v\ \middle|\ \exists\lambda_j\ge0:\
   \sum_j\lambda_j=t,\ \sum_j\lambda_j u_j<v+\one\right\rangle,
\end{equation}
where $v\in\Z_{\ge 0}^n$ and the $\lambda_j$ are real.
The same monomials generate $\tau(\aideal_p^t)$ over $\F_p[\xx]$.
\end{proposition}
For $p\nmid D$, write $\aideal_{f,p}=(\aideal_f)_p=\aideal_{f_p}$.
Since $\langle f_p\rangle\subseteq\aideal_{f,p}$,
\begin{equation}\label{eq:upper}
 \tau(f_p^t)\subseteq\tau(\aideal_{f,p}^{\,t})
      =\J(\aideal_f^t)_p\qquad(t\ge0).
\end{equation}
This inclusion holds for all exponents at once.
\begin{proposition}\label{prop:charzero}
If the support of $f$ is affinely independent, then
\[
 \J(f^t)=\J(\aideal_f^t)\qquad(0\le t<1).
\]
\end{proposition}
\begin{proof}
Work over $\C$. By the general-member formula, a polynomial
$h=\sum_j\eta_j\xx^{u_j}$ with general coefficients satisfies
$\J(h^t)=\J(\aideal_f^t)$ for $t<1$; see \cite[Theorem~2]{How03}.
The torus homomorphism
\[
 (\C^*)^n\times\C^*\longrightarrow(\C^*)^s,
 \qquad(\rho,\delta)\longmapsto(\delta\rho^{u_j})_j
\]
is surjective: its character map sends the standard basis to the
linearly independent vectors $(u_j,1)$.
Choose $(\rho,\delta)$ with $\eta_j=\delta\gamma_j\rho^{u_j}$.
Then $h=\delta\sigma(f)$ for the diagonal automorphism
$\sigma(\xx)=\rho\xx$. This automorphism preserves every monomial ideal,
including $\J(\aideal_f^t)$. Applying $\sigma^{-1}$ gives
the claimed equality over $\C$. Multiplier ideals commute with field
extension in characteristic zero, so faithful flatness gives the equality
over $\Q$.
\end{proof}
%%%%%%%%%%%%%%%%%%%%%%%%%%%%%%%%%%%%%%
%%%%%%%%%%%%%%%%%%%%%%%%%%%%%%%%%%%%%%
\section{Monomials in Frobenius roots}
Throughout this section, let $f$ be as in \eqref{eq:f}
and let $A=\Z[1/D]$ be as in Section~2.

\begin{lemma}\label{lem:closed}
Let $0<\beta\le1$ be rational. If $\xx^v\in\J(\aideal_f^t)$ for all
$0\le t<\beta$, there are rational numbers
$\lambda_j\ge0$ such that
\begin{equation}\label{eq:closed}
 \sum_j\lambda_j=\beta,\qquad \sum_j\lambda_j u_j\le v+\one.
\end{equation}
\end{lemma}
\begin{proof}
Choose a sequence $0\le t_\nu\uparrow\beta$.
For each $\nu$, choose
$(\lambda_1^{(\nu)},\ldots,\lambda_s^{(\nu)})$
satisfying \eqref{eq:monomial} for $t_\nu$.
These vectors lie in the compact set $[0,\beta]^s$.
A convergent subsequence gives a real point of
\[
 \left\{\lambda\in\R_{\ge0}^s:\
  \sum_j\lambda_j=\beta,\ \sum_j\lambda_j u_j\le v+\one\right\}.
\]
This is a nonempty bounded rational polyhedron, so it has a rational
vertex. This vertex gives the required numbers.
\end{proof}
For $m,k_1,\dots,k_s \in \Z_{\ge 0}$ with $m=k_1+\cdots+k_s$, 
$\binom{m}{k_1,\ldots,k_s}$ denotes the multinomial coefficient;
\[
(T_1+\dots+T_s)^m=\sum_{k_1+\dots+k_s =m} \binom{m}{k_1,\ldots,k_s}T_1^{k_1} \cdots T_s^{k_s}. 
\]
The following is the classical multinomial Lucas congruence; see
\cite[Theorem~2.1]{HTWF}, attributed there to Dickson.
\begin{lemma}[Multinomial Lucas congruence]\label{lem:lucas}
Let $\sum_jk_j=m$, and write the base-$p$ expansions as
$k_j=\sum_\nu k_{j,\nu}p^\nu$ and $m=\sum_\nu m_\nu p^\nu$.
Then

\[
 \binom{m}{k_1,\ldots,k_s}
 \equiv\prod_\nu\binom{m_\nu}{k_{1,\nu},\ldots,k_{s,\nu}}\pmod p,
\]
where a factor is zero unless its lower entries sum to its upper entry.
The coefficient is nonzero modulo $p$ if and only if
$\sum_jk_{j,\nu}=m_\nu$ at every digit.
\end{lemma}
\begin{proof}
Expand
$(T_1+\cdots+T_s)^m=\prod_\nu(T_1^{p^\nu}+\cdots+T_s^{p^\nu})^{m_\nu}$
in characteristic $p$ and compare coefficients.
Each $m_\nu<p$, so a digit coefficient is nonzero exactly when its lower
entries sum to $m_\nu$.
\end{proof}

For affinely independent support, let $\Delta$ be the positive greatest
common divisor of the maximal minors of the matrix with columns $(u_j,1)$.

\begin{lemma}[Separation of terms]\label{lem:separation}
Suppose the support is affinely independent and $p\nmid D\Delta$.
For $q=p^e$ with $e\ge1$ and $0\le m<q$ an integer, let $k$ range over the vectors in $\Z_{\ge 0}^s$ with $\sum_j k_j=m$. Then
\[
\langle f_p^m\rangle^{[1/q]}=
\left\langle \xx^{\lfloor\sum_j k_j u_j/q\rfloor}
\ \middle|\ \
\binom{m}{k_1,\ldots,k_s}\not\equiv0\pmod p\right\rangle .
\]
\end{lemma}
\begin{proof}
Some maximal minor is a unit modulo $q$, so the matrix with columns
$(u_j,1)$ is injective modulo $q$. If two multiplicity vectors $k,k'$
of total degree $m$ give the same exponent modulo $q$, this matrix
sends $k-k'$ to zero modulo $q$. Hence $k\equiv k'\pmod q$.
Each entry lies between $0$ and $m<q$, so $k=k'$.
Since each $\gamma_j$ is a unit in $A=\Z[1/D]$ and $p\nmid D$,
all coefficients remain nonzero modulo $p$. 
The multinomial expansion of $f_p^m$ and Lemma~\ref{lem:froot}
now give the formula.
\end{proof}

The proof of the following lemma uses the repeated-digit argument
from \cite[Example~4.2]{MTW}.

\begin{lemma}\label{lem:endpoint}
Suppose the support is affinely independent. Let $v\in\Z_{\ge 0}^n$,
$0<\beta\le1$, and $\lambda_j\ge0$ satisfy \eqref{eq:closed}.
Assume that $p\nmid D\Delta$ and that $(p-1)\lambda_j\in\Z$ for every $j$.
Then for every $q=p^b$ with $b\ge1$,
\[
 \xx^v\in\langle f_p^{\beta(q-1)}\rangle^{[1/q]}.
\]
\end{lemma}
\begin{proof}
%Put $a_j=(p-1)\lambda_j$ and $B=\sum_ja_j=\beta(p-1)$.
%These are nonnegative integers with $B\le p-1$.
For $q=p^b$, 
\[
\begin{aligned} (q-1)\lambda_j&=(p-1)\lambda_j (1+p+\cdots+p^{b-1}),\\
 \beta(q-1)&=\beta(p-1) (1+p+\cdots+p^{b-1}).\end{aligned}
\]
Since $0\le(p-1)\lambda_j\le(p-1)\beta\le p-1$, the first $b$
base-$p$ digits of $(q-1)\lambda_j$ are all $(p-1)\lambda_j$,
and those of $\beta(q-1)$ are all $(p-1)\beta$.
Lemma~\ref{lem:lucas} gives
\[
 \binom{\beta(q-1)}{\lambda_1(q-1),\ldots,\lambda_s(q-1)}
 \equiv\binom{\beta(p-1)}{\lambda_1(p-1),\ldots,\lambda_s(p-1)}^{b}\not\equiv0\pmod p.
\]
Since $\beta(q-1)\le q-1$, it follows that 
\[
\xx^{\left \lfloor \frac1q \sum_j\lambda_j(q-1)u_j\right\rfloor}\in \langle f_p^{\beta(q-1)}\rangle^{[1/q]}
\]
by Lemma~\ref{lem:separation}. 
Its exponent satisfies
\[
 \frac1q\sum_j\lambda_j(q-1)u_j
 =\left(1-\frac1q\right)\sum_j\lambda_ju_j
 \le\left(1-\frac1q\right)(v+\one)<v+\one.
\]
The last inequality is strict because every coordinate of $v+\one$ is
positive. Thus $\left\lfloor\frac1q\sum_j \lambda_j(q-1)u_j\right\rfloor\le v$, so the displayed monomial divides $\xx^v$. This proves the claim.
\end{proof}

\section{Proof of the main theorem}

Adjoin $0$ and $1$ to the jumping numbers of $\J(\aideal_f^t)$ in $(0,1)$,
and list them as $0=\beta_0<\beta_1<\cdots<\beta_r=1$.
This is a finite rational list by the monomial formula.
For each $i$, let $G_i$ be the exponent vectors of the minimal monomial
generators of the constant ideal
\[
 \J(\aideal_f^t)=\J(\aideal_f^{\beta_{i-1}})
                  \qquad(\beta_{i-1}\le t<\beta_i).
\]
For each $v\in G_i$, the same monomial $\xx^v$ belongs to this ideal
as $t$ approaches $\beta_i$ from below. Applying Lemma~\ref{lem:closed}
to this fixed monomial at $\beta_i$, choose rational numbers
$\lambda_j^{i,v}\ge0$ with
\begin{equation}\label{eq:finite-witnesses}
 \sum_j\lambda_j^{i,v}=\beta_i,\qquad
 \sum_j\lambda_j^{i,v}u_j\le v+\one.
\end{equation}
Fix one such vector for each pair $(i,v)$. This finite list depends only on the support.

\begin{proof}[Proof of Theorem~\ref{thm:main}]
Choose $N$ to be the least common multiple of the denominators in
\eqref{eq:finite-witnesses}. For every sufficiently large prime
$p\equiv1\pmod N$, all reductions are defined, the coefficients of $f$
remain nonzero, $p\nmid D\Delta$, and Lemma~\ref{lem:separation} applies. The congruence
gives $(p-1)\lambda_j^{i,v}\in\Z$ for every $i,v,j$.

Fix one constancy interval $[\alpha,\beta)=[\beta_{i-1},\beta_i)$.
Put $q=p^b$ and $\varepsilon_b=1/q$, taking $b$ large enough that
$\alpha\le \beta -\varepsilon_b <\beta$. Since $\beta(q-1)$ is an integer and
$q(\beta -\varepsilon_b)=\beta(q-1)+\beta-1$ with $-1<\beta-1\le0$, we have
$\lceil q (\beta -\varepsilon_b)\rceil=\beta(q-1)$. Thus
\begin{equation}\label{eq:endpoint-chain}
\begin{aligned}
 \J(\aideal_f^\alpha)_p
 &=\J(\aideal_f^{\beta -\varepsilon_b})_p=\J(f^{\beta -\varepsilon_b})_p\\
 &\supseteq\tau(f_p^{\beta -\varepsilon_b})\\
 &\supseteq\langle f_p^{\lceil q(\beta -\varepsilon_b)\rceil}\rangle^{[1/q]}
 =\langle f_p^{\beta(q-1)}\rangle^{[1/q]}.
\end{aligned}
\end{equation}
The first line uses constancy and Proposition~\ref{prop:charzero}.
The inclusions follow from \eqref{eq:upper} and the term $e=b$ in
\eqref{eq:test-union}, respectively.

For each $v\in G_i$, apply Lemma~\ref{lem:endpoint} to the vector
chosen in \eqref{eq:finite-witnesses}. The rightmost ideal contains
every minimal generator of $\J(\aideal_f^\alpha)_p$, so all ideals
in \eqref{eq:endpoint-chain} are equal.

For any $t\in[\alpha,\beta)$, take $b$ large enough that $t\le \beta -\varepsilon_b$.
%Constancy, the equality just proved, monotonicity in the exponent,
%and \eqref{eq:upper} give
Then we conclude that 
\[
 \J(\aideal_f^t)_p=\tau(f_p^{\beta -\varepsilon_b})
 \subseteq\tau(f_p^t)\subseteq \J(f^t)_p=\J(\aideal_f^t)_p.
\]
Hence $\tau(f_p^t)=\J(f^t)_p=\J(\aideal_f^t)_p$
holds for every $t\in[0,1)$.
Equation~\eqref{eq:skoda} now gives the theorem for all exponents.
Dirichlet's theorem \cite{Dav} gives infinitely many primes in this
progression, even after the finite exclusions.
\end{proof}

\begin{remark}
The congruence condition must clear the denominators of the
$\lambda_j^{i,v}$, not only those of the jumping numbers.
Also, a vector chosen at one interior point of an interval need not work
at its endpoint. Lemma~\ref{lem:closed} gives a rational vector satisfying
the weak inequality at $\beta$. Multiplication by $(q-1)/q$ makes it strict.
The same progression works for any nonzero rational coefficients on the
fixed support, apart from their coefficient-dependent exceptional primes.
\end{remark}

\section{An example}

\begin{example}
Let $f=x^4+y^4+xy^2$. The augmented exponent matrix has determinant $4$.
For $i,j\in\Z_{\ge 0}$, let $c(i,j)$ be the optimal value of
\[
 \max_{\lambda_1,\lambda_2,\lambda_3\ge0}
 \{\lambda_1+\lambda_2+\lambda_3\mid
 4\lambda_1+\lambda_3\le i+1,\quad 4\lambda_2+2\lambda_3\le j+1\}.
\]
By Proposition~\ref{prop:monomial}, $x^iy^j\in\J(\aideal_f^t)$
exactly when $0\le t<c(i,j)$. Scaling an optimal vector to any smaller
nonnegative sum makes both inequalities strict. Conversely, if both
inequalities are strict, one can increase $\lambda_1$ slightly,
so the sum is less than the optimum.

The dual problem minimizes $(i+1)a+(j+1)b$ subject to
$4a\ge1$, $4b\ge1$, and $a+2b\ge1$. Its feasible region has vertices
$(1/4,3/8)$ and $(1/2,1/4)$. By linear programming duality,
\[
 c(i,j)=\min\left\{\frac{2(i+1)+3(j+1)}8,\,
                         \frac{2(i+1)+(j+1)}4\right\}.
\]
Hence
\[
 \J(\aideal_f^t)=\begin{cases}
 \Q[x,y],&0\le t<5/8,\\
 \langle x,y\rangle,&5/8\le t<7/8,\\
 \langle x^2,y\rangle,&7/8\le t<1.
 \end{cases}
\]
For each right endpoint $\beta$ and each generator of the ideal on
the interval to its left, the table gives a choice of $\lambda$ in
\eqref{eq:finite-witnesses}, with coordinates ordered as $x^4,y^4,xy^2$.
\[
\begin{array}{c|c|c}
 \beta&\text{generator}&\lambda\\ \hline
 5/8&1&(1/8,0,1/2)\\
 7/8&x&(3/8,0,1/2)\\
 7/8&y&(0,0,7/8)\\
 1&x^2&(1/2,0,1/2)\\
 1&y&(0,0,1)
\end{array}
\]
For every prime $p\equiv1\pmod8$, the determinant is invertible modulo
$p$ and $(p-1)\lambda$ is integral in every row. Thus the displayed
ideals give $\tau(f_p^t)$ after reduction, and Skoda's theorem gives
the comparison for all remaining exponents.

Consider the row for $x$ at $\beta=7/8$. Take $p\equiv1\pmod8$ and
$q=p^b$ with $b\ge1$. For $\lambda=(3/8,0,1/2)$, the term of
$f_p^{7(q-1)/8}$ with multiplicities $(q-1)\lambda$ has monomial part
\[
 (x^4)^{3(q-1)/8}(xy^2)^{(q-1)/2}
 =x^q\bigl(x^{q-2}y^{q-1}\bigr).
\]
Its multinomial coefficient is nonzero modulo $p$ by
Lemma~\ref{lem:lucas}, since the base-$p$ digits add without a carry.
By Lemma~\ref{lem:separation}, this term is alone in its Frobenius basis
component. Its contribution to the Frobenius root is therefore $x$ up
to a unit. Since $\lceil q(7/8-1/q)\rceil=7(q-1)/8$, we obtain
\[
 x\in\langle f_p^{7(q-1)/8}\rangle^{[1/q]}
 \subseteq\tau(f_p^{7/8-1/q}).
\]
For each $0\le t<7/8$, choose $b$ so that $t\le7/8-1/q$.
Since test ideals decrease as the exponent increases,
we have $x\in\tau(f_p^t)$. 
\end{example}

\end{document}